\documentclass{amsart}

\newtheorem{theorem}{Theorem}[section]
\newtheorem*{maintheorem}{Theorem}
\newtheorem{lemma}[theorem]{Lemma}
\newtheorem{proposition}[theorem]{Proposition}
\newtheorem{corollary}[theorem]{Corollary}

\theoremstyle{definition}

\newtheorem{remark}[theorem]{Remark}
\newtheorem*{acknowledgement}{Acknowledgement}

\theoremstyle{remark}

\DeclareFontFamily{U}{wncy}{}
\DeclareFontShape{U}{wncy}{m}{n}{<->wncyr10}{}
\DeclareSymbolFont{mcy}{U}{wncy}{m}{n}
\DeclareMathSymbol{\Sh}{\mathord}{mcy}{"58}

\newcommand\mynote[1]{\marginpar{\ \\ \small \tt #1}}
\newcommand\bel[1]{{\mynote{#1}}\begin{equation}\label{#1}}
\newcommand\mylabel[1]{\label{#1}}

\newcommand{\ZZ}{\mathbb{Z}}

\newcommand{\PP}{\mathbb{P}}
\renewcommand{\AA}{\mathbb{A}}
\newcommand{\GG}{\mathbb{G}}

\newcommand  {\shA}     {\mathcal{A}}

\newcommand  {\shF}     {\mathcal{F}}

\newcommand  {\shHom}   {\mathcal{H}\!\text{\textit{om}}}

\newcommand  {\shN}     {\mathcal{N}}
\newcommand  {\shL}     {\mathcal{L}}

\newcommand  {\edim}    {\operatorname{edim}}

\newcommand  {\Gal}     {\operatorname{Gal}}

\newcommand  {\gr}      {\operatorname{gr}}

\newcommand  {\lra}     {\longrightarrow}

\newcommand  {\maxid}   {\mathfrak{m}}

\renewcommand{\O}       {\mathcal{O}}

\newcommand  {\Pic}     {\operatorname{Pic}}

\newcommand  {\pr}      {\operatorname{pr}}
\newcommand  {\Proj}    {\operatorname{Proj}}

\newcommand  {\quadand} {\quad\text{and}\quad}

\newcommand  {\ra}      {\rightarrow}

\newcommand  {\red}     {{\operatorname{red}}}

\newcommand  {\Sing}    {\operatorname{Sing}}

\newcommand  {\Spec}    {\operatorname{Spec}}

\newcommand  {\Sym}     {\operatorname{Sym}}

\def\mydate{\number\day\space\ifcase\month \or January\or February\or March\or 
April\or May\or June\or July\or
August\or September\or October\or November\or December\fi \space\number\year}

\DeclareFontFamily{U}{wncy}{}
\DeclareFontShape{U}{wncy}{m}{n}{<->wncyr10}{}
\DeclareSymbolFont{mcy}{U}{wncy}{m}{n}
\DeclareMathSymbol{\Sh}{\mathord}{mcy}{"58}

\usepackage{amscd,amssymb,amsmath}

\begin{document}

\title[Fibrations whose geometric fibers are nonreduced]
      {On Fibrations  whose geometric fibers are    nonreduced}


\author[Stefan Schr\"oer]{Stefan Schr\"oer}
\address{Mathematisches Institut, Heinrich-Heine-Universit\"at,
40225 D\"usseldorf, Germany}
\curraddr{}
\email{schroeer@math.uni-duesseldorf.de}

\subjclass{14A15, 14D05, 14J70}

\dedicatory{25 February 2009}

\begin{abstract}
We give a bound on embedding dimensions of geometric generic fibers
in terms of the dimension of the base, for fibrations
in positive characteristic. This generalizes the well-known fact
that for fibrations over curves, the geometric generic fiber is reduced. 
We illustrate our results with Fermat hypersurfaces and genus one curves.
\end{abstract}

\maketitle
\tableofcontents
\renewcommand{\labelenumi}{(\roman{enumi})}

\section*{Introduction}
The goal of this paper is to understand \emph{geometric nonreducedness} for fibrations
in characteristic $p>0$. Roughly speaking, we deal with nonreducedness  that becomes visible only after
purely inseparable base change.

My starting point is the following well-known fact: Let $k$ be an algebraically closed
ground field,   $S$ be a normal algebraic scheme, and $f:S\ra C$ be a fibration
over a curve: then the geometric generic fiber $S_{\bar{\eta}}$ is reduced.
This fact goes back to MacLane  (\cite{MacLane 1939}, Theorem 2); a proof can also be found
in B\u{a}descu's monograph (\cite{Badescu 2001}, Lemma 7.2).
It follows, for example,  that geometric nonreducedness plays no role in the Enriques classification of surfaces.
A natural question comes to mind: Are there generalisations to higher dimensions?
The  main result of this paper is indeed such a generalisation:

\begin{maintheorem}
Let $f:X\ra B$ be a proper morphism with $\O_B=f_*(\O_X)$ between normal   $k$-schemes of finite type.
Then the generic embedding dimension of $X_{\bar{\eta}}$ is smaller the $\dim(B)$.
\end{maintheorem}

Here the \emph{generic embedding dimension} is the embedding dimension of the local ring at the generic
point of $X_{\bar{\eta}}$. 

I expect that the phenomenon of geometric nonreducedness in fibrations  will
play an   role  in the characteristic-$p$-theory of genus-one-fibrations, Albanese maps,
and the minimal model program. For example,
 Koll\'{a}r mentions the threefold $X\subset\PP^2\times\PP^2$ in characteristic two defined by 
the bihomogeneous equation $xu^2+yv^2+zw^2=0$, whose projection onto the first factor defines
a Mori fiber space that has the structure of a geometrically nonreduced conic bundle (\cite{Kollar 1991}, Example 4.12).
The whole theory, however, is largely unexplored to date.

It is actually better to formulate our result in a more general framework, in which an
arbitrary field $K$ takes over the role of the function field of the scheme $B$. In this setting,
the dimension of $B$, which is also the transcendence degree of the function field $\kappa(B)$,
has to be replaced by the \emph{degree of inseparability} of $K$, which is 
the length of any $p$-basis for the extension $K^p\subset K$.

\begin{maintheorem}
Let $X$ be a proper normal $K$-scheme with $K=H^0(X,\O_X)$.
Suppose $X$ is geometrically nonreduced. Then the geometric generic embedding dimension of $X$
is smaller then the degree of inseparability of $K$.
\end{maintheorem}

The key observation of this paper is that  base change with  field extensions $K\subset L$ of degree $p$ cannot
produce nilpotent functions on $X\otimes_KL$. However, the global functions on the \emph{normalisation} of $X\otimes_KL$
may form a field that is larger than $L$.
Our arguments also hinge on  Kraft's beautiful description
of finitely generated field extensions in positive characteristics \cite{Kraft 1970}.

To illustrate our result, we study \emph{$p$-Fermat hypersurfaces} $X\subset\PP^n_K$, which are
defined by a homogeneous equation of the form
$$
\lambda_0U_0^p+\lambda_1U_1^p+\ldots+\lambda_nU_n^p=0.
$$
These are twisted forms of an infinitesimal neighborhood of a hyperplane.
Mori and Saito \cite{Mori; Saito 2003} studied them in the context of fibrations, using  the name   \emph{wild hypersurface bundle},
in connection with the minimal model program.
We attach to such $p$-Fermat hypersurfaces a numerical invariant $0\leq d\leq n$
defined as the $p$-degree of certain field extension   depending on the coefficients
$\lambda_i$ (for details, see Section \ref{fermat hypersurfaces}), and show that this number has a geometric significance:

\begin{maintheorem}
The $p$-Fermat hypersurface $X\subset\PP^n$ is regular if and only if $d=n$.
If $X$ is not regular, then  the singular set $\Sing(X)\subset X$ has codimension $0\leq d\leq n-1$.
\end{maintheorem}

This generalizes parts of  a result of Buchweitz, Eisenbud and Herzog, who studied quadrics
in characteristic two (\cite{Buchweitz; Eisenbud; Herzog 1987}, Theorem 1.1). Our proof uses
different methods and relies on Grothendieck's  
theory of \emph{the generic hyperplane section}.
The result nicely shows that field extension  of degree $p$ may increase the dimension of the singular set
only stepwise.

Then we take a closer look at the   case of\emph{ $p$-Fermat plane curves} $X\subset\PP^2$ that
have an isolated singularity. It turns out that the normalization is a projective line over
a purely inseparable field extension $K\subset L$ of degree $p$. 
Roughly speaking, the curve $X$ arises from this projective line by \emph{thinning out} an infinitesimal
neighborhood of a $L$-rational point. Here some amazing phenomena related to nonuniqueness of 
\emph{coefficient fields in local rings} come into light. 

We then study the question whether curves $C$ arising as such a denormalization
are necessarily $p$-Fermat plane  curve. I could not settle this. However,
we shall see that this is true under the assumption that $p\neq 3$ and  $C$ globally embeds
into a smooth surface.
In this context we touch upon the theory of \emph{abstract multiple curves}, which 
were studied, for example,  by B\u{a}nic\u{a} and Forster \cite{Banica; Forster 1986}, Manolache \cite{Manolache 1994}, and Drezet \cite{Drezet 2007}.
We also make use of Grothendieck's theory of \emph{Brauer--Severi schemes}.

Finally, we study \emph{genus one curves}, that
is, proper curves   $X$ with cohomological invariants $h^0(\O_X)=h^(1\O_X)=1$. They should play an important role in the 
characteristic-$p$-theory of genus one fibrations.  Here we have a result
on the structure of the  Picard scheme:
If  $X$ be a  genus one curve that is regular but geometrically nonreduced, then
 $\Pic^0_X$ is unipotent, that is, a twisted form of the additive group scheme
$\GG_a$.
As a consequence, the reduction of $X\otimes_K\bar{K}$ cannot be 
an elliptic curve.

\begin{acknowledgement}
It is a pleasure to thank Igor Dolgachev: stimulating discussions with him
originated this research.
\end{acknowledgement}

\section{Geometric nonreducedness}
\mylabel{geometric nonreducedness}

Throughout this paper, we fix a prime number $p>0$.
Let $K$ be a field of characteristic $p$,
and let $X$ be a  proper   $K$-scheme so that the canonical inclusion $K\subset H^0(X,\O_X)$
is bijective. The latter can always be achieved, at least if $X$ is connected and reduced, by replacing the
ground field by the finite extension field $H^0(X,\O_X)$.
If $K\subset K'$ is a field extension, the induced proper $K'$-scheme $X'=X\otimes_K K'$
may be less regular than $X$, and in fact nilpotent elements may appear. Such phenomena are the topic of this paper.
We start with    obvious facts:

\begin{proposition}
\mylabel{embedded components}
If $X$ contains no embedded components, then $X'=X\otimes_KK'$ contains no embedded
components. In particular, if $X$ is reduced, then $X'$ is reduced if and only if
$X'$ is generically reduced.
\end{proposition}

\proof
The first statement is contained in \cite{EGA IVb}, Proposition 4.2.7,
and the second is an immediate consequence.
\qed

\medskip
Recall that $X$ is called \emph{geometrically reduced} if $X'=X\otimes_KK'$ is reduced for
all field extension $K\subset K'$.
Actually, it suffices to check reducedness   for a single field extension:

\begin{proposition}
\mylabel{geometrically reduced}
The scheme $X$ is geometrically reduced if and only if $X\otimes_K K^{1/p}$ is reduced.
\end{proposition}

\proof
Necessity is trivial. Now suppose that $X\otimes_K K^{1/p}$ is reduced.
Let $\eta\in X$ be a generic point, and let $F=\O_{X,\eta}$ be the corresponding field of functions.
Then $F\otimes_KK^{1/p}$ is reduced. By MacLane's Criterion (see for example \cite{A 4-7}, Chapter IV, \S15, No.\ 4, Theorem 2),
the field extension $K\subset F$ is separable. In light of Proposition \ref{embedded components},
$X$ must be geometrically reduced.
\qed

\medskip
Purely inspeparable field extension that are smaller then $K^{1/p}$
do not necessarily uncover geometric nonreducedness. In fact, degree $p$ extensions
are incapable of doing so:

\begin{lemma}
\mylabel{remains integral}
Suppose that $K\subset K'$ is purely inseparable of degree $p$.
If the scheme $X$ is normal than the induced scheme $X'=X\otimes_K K'$ is at least integral.
\end{lemma}

\proof
Let $F=\O_{X,\eta}$ be the function field of $X$.
In light of Proposition \ref{embedded components}, it suffices to check that the local Artin ring $F'=F\otimes_KK'$ is a field.
Choose an element $a\in K'$ not contained in $K$. Then $b=a^p$ lies in the subfield $K\subset K'$, and $T^p-b\in K[T]$ is the minimal polynomial of $a\in K'$.
It follows that  $K'=K[T]/(T^p-b)$ and $F'=F[T]/(T^p-b)$. Thus our task is to show that $b\in F$ is not a $p$-th power.
Seeking a contradiction, we assume $b=c^p$ for some $c\in F$. Then for each affine open subset $U=\Spec(A)$ in $X$, the element $c\in F$
lies in the integral closure of $A\subset F$. Since $X$ is normal, $c\in A$, and whence $c\in H^0(X,\O_X)$ according to the sheaf axioms.
By our overall assumption  $K=H^0(X,\O_X)$, and we conclude that $b$ is a $p$-th power in $K$, contradiction.
\qed

\medskip
Suppose our  field extension $K\subset K'$ is so that the induced scheme $X'=X\otimes_KK'$
remains integral; then we consider its  normalization  $Z\ra X'$ and obtain 
a new field $L=H^0(Z ,\O_Z)$, such that we have a sequence  of   field extensions
$$
K\subset K'\subset L.
$$
Here $K'\subset L$ is finite. In some sense and under suitable assumption, this extension $K\subset L$ is also not too big.
Recall that a purely inseparable field extension $K\subset E$ is called \emph{of height} $\leq 1$ if one has $E^{p}\subset K$.

\begin{proposition}
\mylabel{height one}
Suppose that $K\subset K'$ is purely inseparable of height $\leq 1$ so that $X'=X\otimes_KK'$
is integral. Then the
field extension $K\subset L$ is purely inseparable and has height $\leq 1$ as well.
\end{proposition}

\proof
Let $f\in H^0(Z,\O_Z)$. We have to show that $g=f^p$ lies
in the image of $H^0(X,\O_X)$ with respect to the canonical projection $Z\ra X$. 
Let $F $ and $F'$ be the function fields of $X$ and $X'$,
respectively. Then $F'=F\otimes_KK'$, and this is also the function field of $Z$.
The description of $F'$ as tensor product gives $g\in F$. 
To finish the proof, let $x\in X$ be a point of codimension one.
Using that $X$ is normal, it suffices to show that $g\in F$ is contained
in the valuation ring $\O_{X,x}\subset F$. Suppose this is not the case.
Then $1/g\in \maxid_x$. Consider the point $z\in Z$
corresponding to $x\in X$. Then $1/g\in\maxid_z$ since $\maxid_z\cap\O_{X,x}=\maxid_x$.
On the other hand, we have $g\in \O_{Z,z}$ and whence $1\in\maxid_z$, contradiction.
\qed

\medskip
Recall that the \emph{$p$-degree} of an extension $K\subset L$ of height $\leq 1$
is the cardinality of any $p$-basis, confer \cite{A 4-7}, Chapter V, \S 13. The \emph{$p$-degree} $[L:K]_p$ of an arbitrary extension is
defined as the $p$-degree of the height one extension $K(L^p)\subset L$.
The \emph{degree of imperfection} of
a field $K$ is the $p$-degree of $K$ over its prime field. In other words, it is the
$p$-degree of $K^p\subset K$, or equivalently  of $K\subset K^{1/p}$. 

\begin{proposition}
\mylabel{maximal extension} 
There is a   purely inseparable extension
$K\subset K'$ of height $\leq 1$ so that the induced scheme $X'=X\otimes_KK'$ is integral,
and that the field of global functions $L=H^0(Z,\O_Z)$ on the
normalization $Z\ra X'$ is isomorphic to $K^{1/p}$  as an extension of $K$.
\end{proposition}

\proof
This is an application of Zorn's Lemma.
Let $K\subset K_\alpha\subset K^{1/p}$ be the collection of
all intermediate fields with $X\otimes_KK_\alpha$ integral.
We may view this collection as an ordered set, where  the ordering comes
from the inclusion relation. This ordered set is inductive by 
\cite{EGA IVc}, Corollary 8.7.3. Dint of Zorn's Lemma,
we choose a maximal intermediate field $K'=K_\alpha$ and consider  $L=H^0(Z,\O_Z)$. Then  $L\subset K^{1/p}$ by Proposition \ref{height  one}.
Seeking a contradiction, we assume that the latter inclusion is not an equality.
Then there is a purely inseparable extension $K\subset E$
of degree $p$ that is linearly disjoint from $K\subset L$. Consider the composite field $K''=K'\otimes_K E$.
Then $K\subset K''$ is purely inseparable of height $\leq 1$ and strictly larger 
then $K'$, so $X''=X\otimes_KK''$ is generically nonreduced. On the other hand,
$X''$ is birational to 
$$
Z\otimes_{K'}K''=Z\otimes_L (L\otimes_{K'}K'')= Z\otimes_L (L\otimes_K E),
$$
which is integral
by Lemma \ref{remains integral}, contradiction.
\qed

\section{Generic embedding dimension}
\mylabel{generic embedding dimension}

\medskip
Let $K$ be a field of characteristic $p>0$, and $K\subset F$ be
a finitely generated extension field. 
If $K\subset K'$ is purely inseparable, the tensor product
$R=F\otimes_KK'$ is a local Artin ring. We now investigate
its \emph{embedding dimension} $\edim(R)$, that is, the smallest number of generators
for the maximal ideal $\maxid_R$, or equivalently the vector space dimension
of the cotangent space $\maxid_R/\maxid_R^2$ over the residue field
$R/\maxid_R$.

We first    will relate the embedding dimension
$\edim(F\otimes_KK')$, which we regard as  an invariant from 
algebraic geometry, with some  invariants from field theory.
Our analysis hinges on Kraft's beautiful result  on the structure of 
finitely generated field extensions: According to \cite{Kraft 1970},
there is a chain of intermediate fields
\begin{equation}
\label{kraft presentation}
K\subset F_0\subset F_1\subset\ldots\subset F_m=F
\end{equation}
so that $K\subset F_0$ is separable,    each  $F_i\subset F_{i+1}$ is purely inseparable,
and moreover $F_{i+1} $ is generated over $F_i$ by a single element $a_i$ whose minimal polynomial is
of the form $T^{p^{r_i}}-b_i$ with constant term $b_i\in K(F_i^{p^{r_i}})$ and $r_i>0$. We now
exploit the rather special form of the constant terms $b_i$.

\begin{proposition}
\mylabel{integers coincide}
Suppose the extension $K\subset K'$ contains $K^{1/p}$. Then the following integers coincide:
\begin{enumerate}
\item
The embedding dimension of $F\otimes_K K'$.
\item
The number $m$ of purely inseparable field extension in the chain (\ref{kraft presentation}).
\item
The difference between the $p$-degree and the transcendence degree of $K\subset F$.
\end{enumerate}
\end{proposition}

\proof
Recall that the \emph{$p$-degree} of an arbitrary field extension $K\subset F$ is defined
as the $p$-degree of the height $\leq 1$ extension
$K(F^p)\subset F$.
Let $t_1,\ldots,t_n\in F_0$ be a separating transcendence basis over $K$, such that $n$ is
the transcendence degree for $K\subset F$. Clearly, $t_1,\ldots,t_n$ together 
with $a_1,\ldots,a_m$ comprise a $p$-basis for $K(F^p)\subset F$, so the integers in
(ii) and (iii) are indeed the same.

We now check   $\edim(F\otimes_KK')=m$ by induction on $m$.
If $m=0$, then $K\subset F$ is separable, whence $F\otimes_KK'$ is a field.
Suppose now $m\geq 1$, and assume inductively that the local Artin ring $F_{m-1}\otimes_KK'$
has embedding dimension $m-1$.
Write $F=F_{m-1}[T]/(T^{p^r}-b)$ for some $b\in K(F_{m-1}^{p^r})$.
Clearly, $b$ is not a $p$-th power, but it becomes a $p$-th power after tensoring
with $K'$ because $K^{1/p}\subset K'$. Now write
$$
F\otimes_K K'= (F_{m-1}\otimes_KK' )[T]/(T^{p^r}-b\otimes1).
$$
Our claim   then follows from the following Lemma.
\qed

\begin{lemma}
\mylabel{plus one}
Suppose $R$ is a local Noetherian ring in characteristic $p>0$,
and let $A=R[T]/(T^{p^r}-f^p)$ for some integer $r\geq 1$ and some element $f\in R$.
Then we have $\edim(A)=\edim(R)+1$.
\end{lemma}

\proof
We may assume that $\maxid_R^2=0$. Let $\bar{f}$ denote the   class of $f$
in the residue field $k=R/\maxid_R$. Write $\bar{f}=\bar{g}^{p^{m-1}}$ with $\bar{g}\in k$
and $0\leq m-1\leq r$ as large as possible, and choose a representant $g\in R$ for $\bar{g}$.
Then $f^p-g^{p^m}\in\maxid_R^p$, whence $f^p=g^{p^m}$. Set $h=T^{p^{r-m}}-g$, such that
$A=R[T]/(h^{p^m})$. Consider the ring $A'=R[T]/(h)$. This is a \emph{gonflement}
of $R$ in the sense of \cite{AC 8-9}, Chapter IX, Appendix to No.\ 1, and we have $\edim(A')=\edim(R)$ by loc.\  cit.\ Proposition 2.
The element $h\in A$ is nilpotent, and  $A'=A/hA$ by definition. It therefore suffices
to check $h\not\in\maxid_A^2$. Seeking a contradiction, we assume $h\in \maxid_A^2$.
Applying the functor $\otimes_Rk$, we reduce to the case that $R=k$ and $R'=k'$ are fields,
such that $\maxid_A=hA$. It then follows $\maxid_A=\cap_{i\geq 0}\maxid_A^i=0$, such that
the projection $A\ra A'$ is bijective. Taking ranks of these free $R$-modules, we obtain $p^m=1$
and thus $m=0$, which contradicts $0\leq m-1$.
\qed

\medskip
In \cite{Kraft 1970}, the integer in Proposition \ref{integers coincide} is
called the \emph{inseparability} of $K\subset F$. 
We prefer to call it the \emph{geometric embedding dimension} for $K\subset F$, to
avoid confusion with other measure of inseparability and to stress its
geometric meaning. Note that this invariant neither depends on the choice of $K'$
nor on the choice of the particular chain (\ref{kraft presentation}).

We now go back to the setting of algebraic geometry:
Suppose $X$ is an integral proper $K$-scheme with $K=H^0(X,\O_X)$.
Let $F=\O_{X,\eta}$ be its function field. We  call the embedding dimension
of $F\otimes_K K^{1/p}$ the \emph{geometric generic embedding dimension} of $X$.
The main result of the paper is the following:

\begin{theorem}
\mylabel{inseparability degree}
Let $X$ be a proper normal $K$-scheme with $K=H^0(X,\O_X)$.
Suppose that $X$ is not geometrically reduced.
Then the geometric generic embedding dimension of $X$ is smaller
then the degree of imperfection for $K$.
\end{theorem}

\proof
The statement is trivial if the degree of imperfection for $K$ is infinite.
Assume now that this degree of imperfection is finite.
According to Proposition \ref{maximal extension}, there is a purely inseparable extension $K\subset K'$
of height $\leq 1$ so that the induced scheme $X'=X\otimes_KK'$ is integral,
and that its normalization $Z$ has the property that $H^0(Z,\O_Z)=K^{1/p}$.
Note that $K$ cannot be perfect, since $X$ is not geometrically reduced,
and hence $K\neq K'$.
We now use the transitivity properties of tensor products: The scheme
$$
X\otimes_KK^{1/p} = (X\otimes_KK')\otimes_{K'}K^{1/p} = X'\otimes_{K'}K^{1/p}
$$ 
is birational to 
$$
Z\otimes_{K'}K^{1/p}= Z\otimes_{K^{1/p}}(K^{1/p}\otimes_{K'}K^{1/p}).
$$
The    local Artin ring $K^{1/p}\otimes_{K'}K^{1/p}$ has residue field $K^{1/p}$, and its embedding
dimension coincides with the $p$-degree of $K'\subset K^{1/p}$,
by Lemma \ref{embedding dimension} below. This embedding dimension is clearly the geometric generic
embedding dimension of $X$.
The $p$-degree of $K'\subset K^{1/p}$  is strictly smaller then the $p$-degree of $K\subset K^{1/p}$, which coincides with the degree of imperfection for $K$.
\qed

\begin{corollary}
\mylabel{base dimension}
Let $k$ be a perfect field,   $f:X\ra B$ be a proper morphism 
with $\O_B=f_*(\O_X)$ between integral normal algebraic $k$-schemes.
Let  $\eta\in B$ be the generic point and suppose $\dim(B)>0$.
Then the geometric generic embedding dimension of $X_\eta$
is smaller then $\dim(B)$.
\end{corollary}

\proof
The statement is trivial if $X_\eta$ is geometrically reduced.
If $X_\eta$ is not geometrically reduced, then its geometric generic embedding
dimension is smaller than the degree of imperfection of the function field $K=\O_{B,\eta}$.
The latter coincides with $\dim(B)$ by \cite{A 4-7}, Chapter IV, \S 16, No.\ 6,
Corollary 2 to Theorem 4, because the ground field $k$ is perfect.
\qed

\medskip
Let us also restate MacLane's result \cite{MacLane 1939}, Theorem 2 in geometric form:

\begin{corollary}
\mylabel{over curve}
Assumptions as in Corollary \ref{base dimension}.
Suppose additionally that $B$ is a curve.
Then $X_\eta$ is geometrically reduced.
\end{corollary}

\proof
According to the previous Corollary, the geometric generic embedding dimension of $X_\eta$
is zero.
\qed

\medskip
In the proof for Theorem \ref{inseparability degree}, we needed the following fact:

\begin{lemma}
\mylabel{embedding dimension}
Let $K\subset L$ be a finite purely inspeparable extension of height $\leq 1$.
Then the local Artin ring $R=L\otimes_{K}L$ has residue field isomorphic to $L$,
and its embedding dimension equals the $p$-degree of $K\subset L$.
\end{lemma}

\proof
Choose a $p$-basis $a_1,\ldots,a_n\in L$, say with $a_i^p=b_i\in K$.
Then we have $L=K[T_1,\ldots,T_n]/(T_1^p-b_1,\ldots,T_n^p-b_n)$, and consequently
$$
R=L[U_1,\ldots,U_n]/(U_1^p-a_1^p,\ldots,U_n^p-a_n^p).
$$
Clearly, the $U_i-a_i$ are nilpotent and generate the $L$-algebra $R$,
and their residue classes in $\maxid_R/\maxid_R^2$
comprise a vector space basis. Hence $\edim(R)$ equals the $p$-degree of $K\subset L$.
\qed

\section{$p$-Fermat hypersurfaces}
\mylabel{fermat hypersurfaces}

Let $K$ be a ground field of characteristic $p>0$.
In this section, we shall consider  \emph{Fermat hypersurfaces} $X\subset\PP^n$ defined by 
a homogeneous equation of the special form 
$$
\lambda_0 U_0^p+\lambda_1 U_1^p+\ldots+\lambda_n U_n^p=0.
$$
Here we write $\PP^n=\Proj(K[U_0,\ldots,U_n])$, and  $\lambda_0,\ldots,\lambda_n\in K$ are scalars
not all zero. Let us call such subschemes \emph{$p$-Fermat hypersurfaces}.
If the scalars are contained in $K^p\subset K$, then $X$
is the  $(p-1)$-th infinitesimal neighborhood of a hyperplane.
In any case, $X$ is a twisted form of the $(p-1)$-th infinitesimal neighborhood of
a hyperplane, such that $X$ is nowhere smooth. The following is immediate: 

\begin{proposition}
\mylabel{rational point}
We have $X(K)=\emptyset$ if and only if the scalars
$\lambda_0,\ldots,\lambda_n\in K$ are linearly independent over $K^p$.
\end{proposition}

The main goal of this section is to understand the  \emph{singular locus} $\Sing(X)\subset X$,
which comprise the points $x\in X$ whose local ring $\O_{X,x}$ is not regular.
To this end, let
$$
f=\lambda_0 U_0^p+\lambda_1 U_1^p+\ldots +\lambda_n U_n^p\in K[U_0,\ldots,U_n]
$$
 be the homogeneous polynomial defining $X=V_+(f)$, and consider the intermediate field
$$
K^p\subset E\subset K
$$
generated over $K^p$ by all fractions $f(\alpha_0,\ldots,\alpha_n)/f(\beta_0,\ldots,\beta_n)$
with nonzero denominator, and scalars $\alpha_i,\beta_i\in K$. Clearly,
this intermediate field depends only on the closed subscheme $X\subset\PP^n$, and does not change if
one translates $X$ by an automorphism of $\PP^n$. 
A more direct but less invariant description of this intermediate field is as follows:

\begin{proposition}
\mylabel{uncanonical description}
Suppose $\lambda_r\neq 0$. Then the extension $K^p\subset E$
is generated by the fractions $\lambda_i/\lambda_r$, $0\leq i\leq n$.
\end{proposition}

\proof
Let $K^p\subset E'\subset K$ be the intermediate field generated by the $\lambda_i/\lambda_r$.
The inclusion $E'\subset E$ is trivial.
The converse inclusion follows from
$$
f(\alpha_0,\ldots,\alpha_n)/f(\beta_0,\ldots,\beta_n) =
{\sum_{i=0}^n}\alpha_i^p
(\sum_{j=0}^n \frac{\lambda_j}{\lambda_r}\frac{\lambda_r}{\lambda_i}\beta_j^p )^{-1},
$$
where the outer sum is restricted to those indices $0\leq i\leq n$ with $\lambda_i\neq 0$.
\qed

\medskip
We obtain a \emph{numerical invariant} $d=[E:K^p]_p$ for our $p$-Fermat variety $X\subset\PP^n$, 
the $p$-degree of the field extension $K^p\subset E$.
In light of Proposition \ref{uncanonical description}, we have $0\leq d\leq n$.
This numerical invariant has a geometric significance:

\begin{theorem}
\mylabel{codimension singular}
The   scheme $X$ is regular if and only if $d=n$.
If the singular set $\Sing(X)\subset X$ is nonempty, then its codimension  equals the $p$-degree
$d=[E:K^p]_p$.
\end{theorem}

\proof
Let $c\in\left\{0,\ldots,n-1,\infty\right\}$ be the codimension of $\Sing(X)\subset X$.
In the first part of the proof we show that $c\geq d$.
By convention, the case $c=\infty $ means that $\Sing(X)$ is empty.
Without loss of generality we may assume that $\lambda_0=1$,
and that $\lambda_1,\ldots,\lambda_d\in K$ are $p$-linearly independent.
For $d+1\leq j \leq n$, we then may write
$$
\lambda_j = P_j(\lambda_1,\ldots,\lambda_d),
$$
where $P_j(V_1,\ldots,V_d)$ is a polynomial with coefficients in $K^p$ and  of degree $\leq p-1$ in each of the   variables.
Since the scalars $\lambda_1,\ldots,\lambda_d\in K$ are $p$-linearly independent,
there  are derivations $D_i:K\ra K$ with $D_i(\lambda_j)=\delta_{ij}$,
the Kronecker Delta. We may extend them to derivations of degree zero
$D_i:K[X_0,\ldots,X_n]\ra K[X_0,\ldots,X_n]$ sending the variables to zero.
Then the singular locus of $X=V_+(f)$ is contained in the closed subscheme
$S\subset\PP^n$ defined by the vanishing of the homogeneous polynomials
\begin{gather*}
f=U_0^p+\lambda_1 U_1^p+\ldots+\lambda_n U_n^p,\\
D_i(f)=U_i^p + \sum_{j=d+1}^{n}\frac{\partial P_j}{\partial V_i}(\lambda_1,\ldots,\lambda_d)U_j^p,\quad 1\leq i\leq d.
\end{gather*}
Substituting the   relations $D_i(f)$ into the former relation $f=0$,
we infer that $S\subset\PP^n$ is defined by the vanishing of homogeneous polynomials of the form
$$
U_i^p-Q_i(U_{d+1},\ldots,U_n),\quad 0\leq i\leq d.
$$
It suffices to check that $\dim(S)\leq n-1-d$. Suppose this is not the case.
By Krull's Principal Ideal Theorem, the dimensions of $S'=S\cap V_+(U_{d+1},\ldots,U_n)$
is of dimension $\geq 0$. On the other hand, we have
$S'_\red=V_+(U_0,\ldots,U_n)=\emptyset$, contradiction.
This proves that $c\geq d$. In particular, the scheme $X$ is regular provided $d=n$.

To finish the proof we claim that $c=d$ if $d<n$.
We proof the claim by induction on $n$. The case $n=1$ is trivial.
Now suppose $n\geq 2$, and that the claim is true for $n-1$.
Suppose $d<n$, otherwise there is nothing to prove.
Without loss of generality we may assume that $\lambda_n\in K$ is nonzero and $p$-linearly
dependent on $\lambda_0,\ldots, \lambda_{n-1}$.

To proceed we employ Grothendieck's theory of \emph{the generic hyperplane
section} as exposed in the unpublished manuscript \cite{EGA V}, confer
also Jouanoulous's monograph \cite{Jouanolou 1983}. Let $\check{\PP}^n$ be the  scheme of hyperplanes in $\PP^n$,
and $H\subset\PP^n\times\check{\PP}^n$ be the universal hyperplane,
and $\eta\in\check{\PP^n}$ be the generic point.
We then denote by $Y=(X\times\check{\PP}^n)\cap H$ the universal hyperplane section of $X$,
and by $Y_\eta=Y\times_{\check{\PP}^n}\Spec \kappa(\eta)$ the generic hyperplane
section of $X$. Note that $Y_\eta$ is actually a hyperplane section in  $X\otimes_K\kappa(\eta)$,
and that the projection $Y_\eta\ra X$ has highly unusual geometric properties.

We have a closed embedding $Y_\eta\subset H_\eta$, 
and $H_\eta$ is isomorphic to the projective space of dimension $n-1$ over
$K'=\kappa(\eta)$. To describe it explicitely, let
$U_i^*$ be the homogeneous coordinates for $\check{\PP}^n$ dual to the $U_i$.
Then the universal hyperplane $H\subset\PP^n\times\check{\PP}^n$ is defined by
the bihomogeneous equation $\sum U_i\otimes U_i^*=0$, and
the function field of $\check{\PP}^n$ is the subfield
$K'\subset K(U^*_0,\ldots,U_n^*)$ 
generated by the fractions $U^*_i/U_n^*$, $0\leq i\leq n$.
In turn, the generic hyperplane $H_\eta\subset\PP^n_\eta$ is defined
by the homogeneous equation $\sum_{i=0}^n U_i\otimes U_i^*/U_n^*$, such that we have
the additional relation $U_n\otimes 1=-\sum_{i=0}^{n-1} U_i\otimes U_i^*/U_n^*$.
Hence we may write the generic hyperplane $H_\eta=\Proj K'[U'_0,\ldots,U'_{n-1}]$
with $U'_i=U_i\otimes 1$.
In these homogeneous coordinates, the
generic hyperplane section is the $p$-Fermat variety given by $Y_\eta=V_+(f')$ with
$$
f'=\sum_{i=0}^{n-1} (\lambda_i/\lambda_n - (U_i^*/U_n^*)^p)(U'_i)^p.
$$
Now let $c'\in\left\{0,\ldots,n-2,\infty\right\}$ be the codimension of the
singular set $\Sing(Y_\eta)\subset Y_\eta$, and $d'=[E':K'^p]_p$ be the numerical
invariant of the $p$-Fermat variety $Y_\eta\subset H_\eta$. Using that $\lambda_n$ is $p$-linearly dependent
on the $\lambda_1,\ldots,\lambda_{n-1}$, we easily see
that $d=d'$. We now use that the projection $Y_\eta\ra X$ is flat with geometrically regular fibers,
and that its set-theoretical image is the \emph{set of nonclosed points} $x\in X$.
This implies that
$$
c'=\begin{cases}
\infty &\text{if $c=n-1$,}\\
c      &\text{if $c\leq n-2$.}
\end{cases}
$$
We now distinguish two cases:
Suppose first that $d=n-1$, such that $d'=n-1$. We already saw that
this implies that $Y_\eta$ is regular, whence $c'=\infty$, and therefore $c=n-1$.
Now consider the case that $d\leq n-2$. Then $d'\leq (n-1)-1$, and
the induction hypothesis implies $c'=d'$, and finally $d=d'=c'=c$.
Summing up, we have $c=d$ in both cases.
\qed

\section{Singularities on $p$-Fermat plane curves}
\mylabel{fermat curves}

Keeping the assumptions of the previous section, we now
study the singularities lying on \emph{$p$-Fermat plane curves}  $X\subset\PP^2_K$,
say defined by 
\begin{equation}
\label{fermat equation}
\lambda_0U_0^p+\lambda_1U_1^p+\lambda_2U_2^p=0
\end{equation}
with $\lambda_0,\lambda_1,\lambda_2\in K$.
First note that $X$ is a geometrically irreducible projective curve,
and has $h^0(\O_X)=1$ and $h^1(\O_X)=(p-1)(p-2)/2$.
Let $0\leq d\leq 2$ be the numerical invariant of the $p$-Fermat variety $X$, as defined
in the previous section. The curve $X$ is regular if $d=2$,   nonreduced if $d=0$,
and has isolated singularities if $d=1$.

With respect to singularities, only the case $d=1$ is of interest.
Throughout, we suppose that $d=1$. Without loss of generality,
we may then assume that $\lambda_2=1$, that $\lambda=\lambda_0$ does not lie in $K^p$,
and that $\lambda_1=P(\lambda)$ is a polynomial of degree $<p$ in $\lambda$.
So $X\subset\PP^2_K$ is defined by the homogeneous equation
$$
\lambda U_0^p+P(\lambda)U_1^p+U_2^p=0.
$$
Consider the field extension $L=K(\lambda^{1/p})$, and let $T_0,T_1$ be two indeterminates.
The map
\begin{gather*}
K[U_0,U_1,U_2]\lra L[T_0,T_1],\\
U_0\longmapsto  T_0 ,\quad
U_1\longmapsto  T_1,\quad
U_2\longmapsto  -\lambda^{1/p}T_0-P(\lambda^{1/p})T_1
\end{gather*}
defines a morphism $\PP^1_L\ra\PP^2_K$ factoring over $X\subset\PP^2_K$.

\begin{proposition}
\mylabel{normalization}
The induced morphism $\nu:\PP^1_L\ra X$ is the normalization of $X$.
\end{proposition}

\proof
It suffices to check that the morphism
$\nu:\PP^1_L\ra X$ has degree one. The intersection $X\cap V_+(U_0)$
is Cartier divisor of length $p$.
Its preimage on $\PP^1_L$ is given by $V_+(T_0)$, which also has length $p$ over $K$.
Whence the degree in question is one.
\qed

\medskip
It follows that the field extension $K\subset L$ is nothing but the field of global 
section on the normalization of $X$. Therefore, it depends only on the curve $X$,
and not on the chosen Fermat equation (\ref{fermat equation}).

\begin{proposition}
\mylabel{singular point}
The singular locus $\Sing(X)$ consists of precisely one point
$a_0\in X$. The corresponding point $a\in\PP^1_L$ on the normalization
has residue field $\kappa(a)=L$.
\end{proposition}

\proof
Taking the derivative with respect to $\lambda$, we see that the singular locus $\Sing(X)$ is contained in the closed subscheme defined by
$$
\lambda U_0^p+P(\lambda)U_1^p+U_2^p=0\quadand
U_0^p+P'(\lambda)U_1^p=0.
$$
Substituting the latter into the former, we see that there is only one singular point $a_0\in X$.
The preimage on $\PP^1_L$ is defined by $T_0^p+P'(\lambda)T_1^p$,
which clearly defines an $L$-rational point $a\in\PP^1_L$. 
\qed

\begin{remark}
\mylabel{homogeneous coordinates}
The homogeneous coordinates for the preimage $a\in\PP^1_L$ of the singular point $a_0\in X$
are $(-P'(\lambda^{1/p}):1)$
\end{remark}

\medskip
Let $\mathfrak{c}\subset\nu_*(\O_{\PP^1_L})$ be the conductor ideal for the finite birational map $\nu:\PP^1_L\ra X$,
which is the largest $\nu_*(\O_{\PP^1_L})$-ideal contained in the subring $\O_X\subset\nu_*(\O_{\PP^1_L})$.
The conductor ideal defines closed subschemes $A\subset\PP^1_L$ and $A_0\subset X$,
such that we have a commutative diagram
\begin{equation}
\label{conductor square}
\begin{CD}
A @>>>\PP^1_L\\
@VVV @VV\nu V\\
A_0 @>>> X,
\end{CD}
\end{equation}
which is cartesian and cocartesian.
By abuse of notation, we write   $\O_{A_0} $ and $\O_A$ for the local Artin rings defining
the conductor schemes $A_0$ and $A$. Then we have $\O_A=L[u]/(u^l)$ for some integer $l\geq 1$, where $u\in\O_{\PP^1,a}$ denotes a uniformizer,
and the subring $\O_{A_0}\subset\O_A$ is a $K$-subalgebra.

\begin{proposition}
\mylabel{properties subalgebras}
We have $\O_A=L[u]/(u^{p-1})$, and $\O_{A_0}\subset\O_A$ is a $K$-subalgebra generated by two
elements so that $L\not\subset\O_{A_0}$ and $\dim_K(\O_{A'})=p(p-1)/2$.
\end{proposition}

\proof
The conductor square (\ref{conductor square}) yields an exact sequence of coherent sheaves
\begin{equation*}
\label{conductor sequence}
0\lra\O_X\lra\O_{A_0}\oplus\nu_*(\O_{\PP^1_L})\lra\nu_*(\O_A)\lra 0,
\end{equation*}
which in turn gives a long exact sequence of $K$-vector spaces
\begin{equation}
\label{exact sequence}
0\lra H^0(X,\O_X)\lra \O_{A_0}\oplus L\lra \O_A\lra H^1(X,\O_X)\lra 0.
\end{equation}
Using that $h^0(\O_X)=1$ we infer that $\O_{A_0}\cap L=K$.
The latter is equivalent to $L\not\subset \O_{A_0}$, because the field
extension $K\subset L$ has prime degree.
Being a complete intersection, the curve $X$ is Gorenstein.
According to \cite{Serre 1975}, Chapter IV,  \S 3.11, Proposition 7, this implies $\dim_K(\O_{A_0})=\dim_K(\O_{A})/2$.
Now set $r=\dim_K(\O_{A})$.
In light of  $h^1(\O_X)=(p-1)(p-2)/2$, it follows $r=(p-1)(p-2)/2 + r/2 + p -1$, whence $r=p(p-1)$.
Finally, since the scheme $\Spec(\O_{A'})$ contains only one point and  embeds into $\PP^2_K$, 
it   embeds even into $\AA^2_K=\Spec K[U_1,U_2]$, whence the $K$-subalgebra $\O_{A_0}\subset \O_{A}$
is generated by two elements.
\qed

\medskip
The situation is very simple in characteristic two:

\begin{corollary}
\mylabel{conductor 2}
If $p=2$, then $\O_A=L$ and $\O_{A_0}=K$.
\end{corollary}

There is more to say for  odd primes. Here we have to distinguish two cases,
according to the residue field  of the singularity $a_0\in X$, which is either
$K$ or $L$:

\begin{proposition}
\mylabel{conductor L}
Suppose $p\geq 3$ and $\kappa(a_0)=\O_{A'}/\maxid_{A'}$ equals $L$. 
Then we have $\O_{A'}=K[\mu+f,g]$ for some $\mu\in L\smallsetminus K$ and 
$f\in\maxid_{A}\smallsetminus \maxid_{A}^2$ and $g\in\maxid_A^2\smallsetminus\maxid_R^3$.
\end{proposition}

\proof
According to \cite{Matsumura 1980}, Theorem 60, the projection $\O_{A'}\ra L$ onto the residue field  admits a section.
In other words, it is possible to embed   $L=\O_{A'}/\maxid_{A'}$
as a \emph{coefficient field}  $L'\subset\O_{A'}$;
but note that such coefficient fields are not unique.
By Proposition \ref{properties subalgebras}, the  $K$-algebra $\O_{A'}$ is generated by two elements, say $h,g\in \O_{A_0}$.
In case $h,g\in\maxid_{A_0}\cup K$ the residue field would be $K$, contradiction.
Without loss of generality we may assume that $h\in\O_{A_0}$ generates a coefficient
field $L'\subset\O_{A_0}$. 
Let $\mu\in L$ be the image of $h$ in the residue field.
Inside $\O_A=L[u]/(u^{p-1})$, we have $h=\mu+f$ for some 
$f\in\maxid_A$. 
To continue, write $g=\epsilon u^l$ for some
unit $\epsilon\in \O_A$ and some integer $0\leq l\leq p-1$.
Adding some polynomial in $h$ to $g$, we may assume $l\geq 1$.
We now check that $l=2$. Clearly, $\O_{A_0}=L'[g]/(g^{d})$,
where $d=\lceil(p-1)/l\rceil$, such that 
$$
p(p-1)/2=\dim_K(\O_{A_0})=p \lceil(p-1)/l\rceil.
$$
This equation implies  $(p-1)/2\geq (p-1)/l> (p-1)/2-1$,
which easily gives us $l=2$.

It remains to verify $f\not\in\maxid_A^2$. 
One way of seeing this is to use that the local Artin ring $\O_A=\O_{\PP^1_L}/\mathfrak{c}\simeq L[u]/(u^{p-1})$
contains a \emph{canonical} coefficient field, namely  the image of $H^0(\PP^1_L,\O_{\PP^1_L})$.
This gives a   splitting $\O_A=L\oplus\maxid_A$ of $L$-vector spaces,
and whence a projection $\pr:\O_A\ra\maxid_A$. 
Note that this projection does not depend on any choices of the uniformizer $u$.
Seeking a contradiction, we now assume that $f\in\maxid_A^2$.
Then $\pr(\O_{A_0})\subset\maxid_A^2$.
On the other hand, the discussion at the beginning of
Section \ref{fermat curves} tell us  that $C\subset\PP^2_K$
is defined by a homogeneous equation of
the form $\lambda U_0^p+P(\lambda)U_1^p+U_2^p=0$.
The preimage  $a\in \PP^1_L$ of the singular point $a_0\in C$   has
homogeneous coordinates  $(-P'(\lambda^{1/p}):1)$, according
to Remark \ref{homogeneous coordinates}.
The normalization is described in Proposition \ref{normalization},
and sends $U_0/U_1$ to the element $T_0/T_1=u-P'(\lambda^{1/p})$,
where $u=T_0/T_1+P'(\lambda^{1/p})$. It follows that $u\in\maxid_A$ is contained in
$\pr(\O_{A_0})\subset\maxid_A$, contradiction.
\qed

\begin{proposition}
\mylabel{conductor K}
Suppose $p\geq 3$ and $\kappa(a_0)=\O_{A_0}/\maxid_{A_0}$ equals $K$. Then we have
$\O_{A_0}=K[v,w]$ for some  $v,w\in\maxid_A$ so that their classes mod $\maxid_A^2$
are $K$-linearly independent.
\end{proposition}

\proof
 Write $\O_{A_0}=K[v,w]$ for some $v,w\in \maxid_A$.
Then the monomials $v^iw^j$ with $i+j\geq p-1$ vanish.
Since we have $\dim_K(\O_{A_0})=p(p-1)/2$, the monomials $v^iw^j\in R'$ with $i+j\leq p-2$
form a $K$-basis. Whence $v^{p-2},w^{p-2}\neq 0$, such that $v,w\not\in\maxid_A^2$.
Seeking a contradiction, we now assume that $v\equiv \alpha w$ modulo $\maxid_A^2$
for some $\alpha\in K$. Then $v^{p-2}=\alpha^{p-2}w^{p-2}$, contradicting 
linear independence.
\qed

\section{Abstract multiple curves}
\mylabel{multiple curves}

We keep the notation from the preceding section, but
now now reverse the situation:
Fix an $L$-rational point $a\in\PP^1_L$.
Let $A\subset\PP^1_E$ be the $(p-2)$-th infinitesimal neighborhood of $a$, and
write $\O_A=L[u]/(u^{p-1})$. Now choose a $K$-subalgebra $\O_{A_0}\subset\O_A$
and consider the resulting morphism $A=\Spec(\O_A)\ra\Spec(\O_{A_0})=A_0$.
The pushout square 
$$
\begin{CD}
A @>>> \PP^1_L\\
@VVV @VVV\\
A_0 @>>> C
\end{CD}
$$
defines a proper integral curve $C$, with normalization $\nu:\PP^1_L\ra C$
and containing a unique singular point $a_0\in C$.

We call the subalgebra $\O_{A_0}\subset\O_A$ \emph{admissible} if
it takes the form described in  Corollary \ref{conductor 2} or Proposition \ref{conductor L} or Proposition \ref{conductor K}.
From now on we assume that our subalgebra is admissible, and   ask whether or not
the resulting curve $C$ is embeddable into $\PP^2_K$ as a $p$-Fermat plane curve. The cohomology groups indeed have the right dimensions:

\begin{proposition}
\mylabel{cohomology dimensions}
We have $h^0(\O_C)=1$ and $h^1(\O_C)=(p-1)(p-2)/2$.
\end{proposition}

\proof
It is easy to see that an admissible subalgebra $\O_{A_0}\subset\O_A$ satisfies the conclusion
of Proposition \ref{properties subalgebras}, that is,
$\dim_K(\O_{A_0})=p(p-1)/2$ and $L\not\subset\O_{A_0}$.
The the statement in question  follows from the
exact sequence
$$
0\lra H^0(C,\O_C)\lra \O_{A_0}\oplus L\lra \O_A\lra H^1(C,\O_C)\lra 0,
$$
as in the proof for Proposition \ref{properties subalgebras}.
\qed

\medskip
To proceed, we have to  study the behavior of $C$ under base change.
If $K\subset K'$ is a field extension, then the induced curve $C'=C\otimes_KK'$
sits inside the cartesian and cocartesian diagram
\begin{equation}
\label{base change diagram}
\begin{CD}
A' @>>> \PP^1_{L\otimes_KK'}\\
@VVV @VVV\\
A_0' @>>> C',
\end{CD}
\end{equation}
with $A'=A\otimes_KK'$ and $A_0'=A_0\otimes_KK'$.
In particular, $C'$ is nonreduced if and only if the extension field $K\subset K'$ 
contains $L$. The next result implies that under suitable assumptions,
$C'$ locally looks like a Cartier divisor inside a regular surface:

\begin{proposition}
Suppose $L\subset K'$. Then every closed point on $C'$ has embedding dimension two.
\end{proposition}

\proof
The assertion is clear outside the closed point $a_0'\in C'$ corresponding to the singularity $a_0\in C$.
To understand the embedding dimension of $\O_{C',a_0'}$, let $M,G$ be two indeterminates, and consider
the $K$-algebra
$$
R=K[M][[G]],
$$
which is a regular   ring of dimension two that is a formally smooth $K$-algebra.
Therefore, it suffices to construct a surjection $R\ra\O_{C,a_0}^\wedge$.
We do this   for the case that $\kappa(a_0)=L$, that is, $\O_{A'}=K[\mu+f,g]$
as in Proposition \ref{conductor L} (the other cases are similar,   actually simpler).
Choose lifts $\tilde{f},\tilde{g}\in\O_{C,a_0}$ for $f,g\in\O_{A'}$, and define
$$
h:R\lra\O_{C,a_0}^\wedge,\quad M\longmapsto \mu+\tilde{f},\quad G\longmapsto\tilde{g}.
$$
Note that we have to work with $R$ rather than its formal completion $K[[M,G]]$,
because the image $\lambda+\tilde{f}$ does not lie in the maximal ideal.
Obviously, the composite map $R\ra\O^\wedge_{C,a_0}\ra\O_{A'}$ is surjective.
Now let $\lambda u^i\in\O_{C,a_0}^\wedge\subset L[[u]]$ be a monomial with $\lambda\in L$ and $i\geq p-1$.
By completeness, it suffices to check that this monomial  lies in $h(R)$ modulo $u^{i+1}$.
Consider first the case that $i$ is even.
Write $g^{i/2}=\lambda'u^i$ modulo $u^{i+1}$ for some nonzero $\lambda'\in L$, and write $\lambda/\lambda' =P(\mu)$
as a polynomial of degree $<p$ with coefficients in $K$ in terms of the generator $\mu\in L$.
Then $\lambda u^i = P(\lambda +f)\cdot g^{i/2}$ modulo $u^{i+1}$.
Finally suppose $i$ is odd. Write $i=p+j$ for some even $j\geq 0$.
Then $(\mu+f)^p=\mu^p+\alpha u^p$ modulo $u^{p+1}$ with nonzero $\mu^p,\alpha\in K$.
Moreover, $g^{j/2}=\lambda' u^j$ modulo $u^{j+1}$ for some $\lambda'\in L$. As above, we find some polynomial $P$ of degree $<p$
with $P(\lambda)=\lambda/\lambda'$. Then 
$$
\lambda u^i=P(\lambda+f) \frac {(\mu+f)^p -\mu^p}{\alpha} g^{j/2} \mod u^{i+1}.
$$
Using that $R$ is $G$-adically complete and that $\O_{C,a_0}^\wedge$ is complete,
we infer that $h:R\ra\O_{C,a_0}^\wedge$ is surjective.
\qed

\begin{proposition}
Suppose   $L\subset K'$. Then $C'_\red=\PP^1_{K'}$.
\end{proposition}

\proof
It suffices to treat the case $K'=L$.
The diagram (\ref{base change diagram}) yields a birational morphism
$$
\nu:\PP^1_{K'}=(\PP^1_{L\otimes_KK'})_\red\lra C'_\red.
$$
Let $a_0'\in C'$ and $a'\in\PP^1_{K'}$ be the points corresponding
to the singularity $a_0\in C$. 
It suffices to check that the fiber $\nu^{-1}(a_0')\subset\PP^1_{K'}$
is the reduced scheme given by $a'$.
In other words, the maximal ideal $\maxid_{A'}\subset\O_{A'}$ is generated
by the maximal ideal $\maxid_{A_0}$ and the nilradical of $\O_{\PP^1_{L\otimes_KK'}}$.
Let us do this in the case the residue field of $a_0\in C$ is $L$,
that is, $\O_{A_0}=K[\mu+f,g]$ as in Proposition 
\ref{conductor L} (the other two cases being similar).
Then $\O_{A_0'}$ is generated over $L$ by $\mu\otimes1+f\otimes1$ and $g\otimes 1$,
whence also by $\mu\otimes1-1\otimes\mu +f\otimes 1$ and $g\otimes 1$.
But the nilradical of $\PP^1_{L\otimes_KK'}$ is generated by $\mu\otimes 1-1\otimes\mu$.
By assumption, $f\in\O_{A}$ generates the maximal ideal,
and this implies that  $f\otimes1\in\O_{A'} $ generates the maximal ideal modulo
$\mu\otimes1-1\otimes\mu$, whence the claim.
\qed

\medskip
If $L\subset K'$, then $C'$ is a nonreduced curve with smooth reduction and multiplicity $p$, and for each closed point $y\in C'$,
the complete local ring $\O_{C',y}^\wedge$ is a quotient of a formal power series ring in two variable over $\kappa(y)$.
Loosing speaking, one may say that $C'$ \emph{locally embeds into regular surfaces}.
Such curves were studied, for example,  by B\u{a}nic\u{a} and Forster \cite{Banica; Forster 1986}, Manolache \cite{Manolache 1994}, and Drezet \cite{Drezet 2007},
and are called \emph{abstract multiple curves}.
Let $\shN\subset\O_{C'}$ be the nilradical. Then the ideal powers
$$
\O_{C'}=\shN^0\supset\shN\supset\shN^2\supset\ldots\supset\shN^p=0
$$ 
define a filtration on the structure sheaf $\O_{C'}$.
The graded piece $\shL=\shN/\shN^2$ is an coherent sheaf on $C'_\red=\O_{\PP^1_{K'}}$,
and we obtain an algebra map $\Sym(\shL)\ra\gr(\O_{C'})$.
A computation in the complete local rings shows that this map is surjective, with
kernel the ideal generated by $\Sym^p(\shL)$, and that $\shL$ is invertible. Set $d=\deg(\shL)$. Using
$$
1-(p-1)(p-2)/2=\chi(\O_{C'})=\chi(\gr(\O_{C'})=\sum_{j=0}^{p-1} (jd+1) = d p(p-1)/2 + p,
$$
we infer $\deg(\shL)=-1$. This shows:

\begin{proposition}
\mylabel{graded algebra}
If we have $L\subset K'$, then the associated graded algebra is given by $\gr(\O_{C'})=\bigoplus_{i=0}^{p-1}\O_{\PP^1_{K'}}(-i)$,
with the obvious multiplication law.
\end{proposition}

We shall say that $C'$   \emph{globally embeds into a smooth surface} if there is a smooth proper
connected $K'$-surface $S$   into which $C'$ embeds. 
We remark in passing  that the surface $S$ is then geometrically connected, or equivalently $K'=H^0(S,\O_S)$,
because we have $K'\subset H^0(S,\O_S)\ra H^0(C',\O_{C'})=K'$.

\begin{theorem}
\mylabel{surface plane}
Suppose  that $C'$ embeds globally into a smooth  surface and that $p\neq 3$. 
Then $C'$ embeds as  a $p$-Fermat plane curve into $\PP^2_{K'}$.
\end{theorem}

\proof
We first consider the case that $L\subset K'$, for which the assumption $p\neq 3$ plays no role.
Choose an embedding $C'\subset S$ into a smooth proper connected surface.
Then $D=C'_\red$ is isomorphic to a projective line. By Proposition \ref{graded algebra}, its selfintersection
 inside $S$ is the number $D^2=1$. Set $\shL=\O_S(D)$. The exact sequence of sheaves
$0\ra\O_S\ra\shL\ra\shL_D\ra 0$ gives an exact sequence
$$
0\ra K'\lra H^0(S,\shL)\lra H^0(D,\shL_D) \lra H^1(S,\O_S).
$$ 
Suppose for the moment that $H^1(S,\O_S)=0$. Then $\shL$ is globally generated and has $h^0(\shL)=3$.
The resulting morphism $r:S\ra\PP^2_{K'}$ is surjective, because $(\shL\cdot\shL)\neq 0$,
and its degree equals $\deg(r)=(\shL\cdot\shL)=1$, whence $r$ is birational.
Moreover, the induced morphism $r:D\ra\PP^2_{K'}$ is a closed embedding,
and $r(D)\subset\PP^2_{K'}$ is a line. Using $D^2=1=r(D)^2$, we conclude that the exceptional curves
for $r$ are disjoint from $D$. It follows that $C'$ embeds into $\PP^2_{K'}$.
Obviously, it becomes a $p$-Fermat plane curve, because $D=C'_\red$ becomes a line.

We now check that indeed $H^1(S,\O_S)=0$. For this we may assume that $K'$ is algebraically closed.
Now $K_S\cdot D=-3$, whence $K_S$ is not numerically effective. By the Enriques classification of surfaces,
$S$ is either the projective plane or ruled. If there is a ruling $f:S\ra B$, then
$D$ is not contained in a fiber, because $D^2>0$. Whence $D\ra B$ is dominant, and it follows
that $B$ is a projective line. In any case, $H^1(S,\O_S)=0$.

Finally, we have to treat the case that $K\subset K'$ is linearly disjoint form $L$.
Tensoring with $K'$, we easily reduce to the case $K'=K$.
Choose   separable and algebraic closures $K\subset K^s\subset\bar{K}$
and an embedding $C\subset S$
into a proper smooth connected surface $S$. 
Following \cite{GB III}, Section 5, we write $\Pic(S/K)=\Pic_{S/K}(K)$ for the group
of rational points on the Picard scheme.
According to the  previous paragraph   $S_{\bar{K}}$ is a rational surface, whence the abelian group $\Pic(S_{\bar{K}})$ is
a free   of finite rank, and the scheme $\Pic_{S/K}$ is \'etale at each point.
It follows   that the canonical map $\Pic(S_{K^s}/K^s)\ra\Pic(S_{\bar{S}}/\bar{S})$ is bijective.
Now let $[D]\in\Pic(S_{K^s})$ be the class
of $D=(C_{\bar{K}})_\red$.  This class is necessarily invariant under
the action of the Galois group $\Gal(K^s/K)$, because $pD=C_{\bar{K}}$ comes
from a curve on $S$ and the abelian group $\Pic(S_{\bar{K}})$ is torsion free. We conclude that the class $[D]\in\Pic(S/K)$ exists, although it does not
come from an invertible sheaf on $S$. However, it gives rise to a 2-dimensional Brauer--Severi scheme $B$,
which comes along with a morphism $r:S\ra B$ that induces the morphism
$r_{\bar{K}}:S_{\bar{K}}\ra\PP^2_{\bar{K}}$ defined above.
The upshot is that there is an embedding $C\subset B$,
and that the class of $C$ inside $\Pic(B/K)=\ZZ$ equals $p$.
On the other hand, the class of the dualizing sheaf $\omega_B$ equals $-3$.
Since $p\neq 3$ by assumption, we must have $\Pic(B)=\ZZ$,
and consequently our Brauer--Severi scheme is $B\simeq\PP^2_K$.
It is then easy to check that $B\subset\PP^2_K$ is indeed a $p$-Feramt plane curve.
\qed

%
%
%

\section{Genus one curves}
\mylabel{genus one}

The goal of this section is to study geometric nonreducedness for genus one curves.
Throughout, $K$ denotes a ground field of arbitrary characteristic $p>0$.
A \emph{genus one curve} is a proper geometrically irreducible curve $X$  over $K$ with $h^0(\O_X)=h^1(\O_X)=1$.
Clearly, this notion is stable under field extensions $K\subset K'$.
Since $h^0(\O_X)=1$, the curve $X$ contains no embedded component, such that
the dualizing sheaf $\omega_X$ exists.

\begin{proposition}
\mylabel{dualizing sheaf}
Let $X$ be a  reduced genus one curve. Then $X$ is Gorenstein, and $\omega_X\simeq\O_X$.
\end{proposition}

\proof
We have $h^0(\omega_X)=h^1(\omega_X)=1$, hence $\omega_X$ admits a nonzero section $s$.
The map $s:\O_X\ra\omega_X$ is injective, because $X$ is reduced.
Hence we have a short exact sequence
$$
0\lra\O_X\stackrel{s}{\lra}\omega_X\lra\shF\lra 0
$$
for some torsion sheaf $\shF$. Using $h^1(\shF)=0$ and $\chi(\shF)=\chi(\omega_X)-\chi(\O_X)=0$,
we conclude $\shF=0$, and the result follows.
\qed

\medskip
Since $H^2(X,\O_X)=0$, the Picard scheme $\Pic_X$ is smooth and 1-dimensional, so the connected component
of the origin $\Pic_X^0$ is either an elliptic curve, a twisted form of $\GG_m$, or a twisted form of $\GG_a$.
In other words, it is proper, of multiplicative type, or unipotent.
A natural question is whether all three possibilities occur in genus one curves that are regular
but geometrically nonreduced. It turns out that this is not the case:

\begin{theorem}
\mylabel{picard unipotent}
Let  $X$ be a  genus one curve that is regular but geometrically nonreduced.
Then $\Pic^0_X$ is unipotent.
\end{theorem}

\proof
Seeking a contradiction, we assume that the Picard scheme is not unipotent.
Choose an algebraic closure $K\subset\bar{K}$, and set $Y=X\otimes_K\bar{K}$.
Proposition \ref{dualizing sheaf} implies that $Y$ is a genus one curve with $\omega_Y=\O_Y$.
Let $\shN\subset\O_Y$ be the nilradical, which defines the closed subscheme
$Y_\red\subset Y$. We first check that $Y_\red$ is also a genus one curve with $\omega_{Y_\red}=\O_{Y_\red}$.
The short exact sequence $0\ra\shN\ra\O_Y\ra\O_{Y_\red}\ra 0$ yields a long exact sequence
$$
H^1(Y,\shN)\lra H^1(Y,\O_Y)\lra H^1(Y_\red,\O_{Y_\red})\lra 0.
$$
Since the Picard scheme $\Pic^0_Y$ contains no unipotent subgroup scheme,   
the restriction mapping $H^1(Y,\O_Y)\ra H^1(Y_\red,\O_{Y_\red})$ is injective ,
which follows from  Proposition \cite{Bosch; Luetkebohmert; Raynaud 1990}, Section 9.2, Proposition 5. Whence $h^1(\O_{Y_\red})=1$.
Furthermore, we have $h^0(\O_{Y_\red})=1$ since $Y_\red$ is integral and $\bar{K}$ is algebraically closed. 
Consequently, $Y_\red$ is a genus one curve,
and Proposition \ref{dualizing sheaf} tells us that $\omega_{Y_\red}=\O_{Y_\red}$.

Relative duality for the inclusion morphism $Y_\red\ra Y$ yields the formula
$$
\O_{Y_\red}=\omega_{Y_\red} =\shHom(\O_{Y_\red},\omega_Y)=\shHom(\O_{Y_\red},\O_Y).
$$
The term on the right is nothing but the annulator ideal $\shA\subset\O_Y$ of the nilradical $\shN\subset\O_Y$,
such that $\shA=\O_{Y_\red}$ as $\O_Y$-modules, and in particular $h^0(\shA)=1$.
To finish the proof, consider the closed subscheme $Y'\subset Y$ defined by $\shA\subset\O_Y$.
By assumption, $\shN\neq 0$, such that $\shA\neq \O_Y$, and therefore $Y'\neq\emptyset$.
We have an exact sequence
$$
0\lra H^0(Y,\shA)\lra H^0(Y,\O_Y)\lra H^0(Y',\O_{Y'}),
$$
and consequently $h^0(\O_Y)\geq 2$, contradiction.
\qed

\begin{corollary}
Let  $X$ be a  genus one curve that is regular but geometrically nonreduced.
Then the reduction of the induced curve $\bar{X}$ over the algebraic closure $K\subset \bar{K}$
is isomorphic to the projective line or the rational cuspidal curve.
\end{corollary}

\begin{proof}
Set $C\subset\bar{X}$ be the reduction. Clearly, $h^1(\O_C)\leq 1$.
We have $C=\PP^1_{\bar{K}}$ if $h^1(\O_C)=0$. Now assume that $h^1(\O_C)=1$.
Then $\Pic^0_C$ is unipotent by Theorem \ref{picard unipotent}.
According to \cite{Bosch; Luetkebohmert; Raynaud 1990}, Section 9.3, Corollary 12, this implies that the normalization map $\nu:C'\ra C$
is a homeomorphism, which means that $C$ is the rational cuspidal curve.
\end{proof}


\end{document}